\documentclass[12pt]{amsart}
\usepackage{amssymb}
\usepackage{url}
\usepackage{comment}

\theoremstyle{plain}
\newtheorem{thm}{Theorem}[section]

\newtheorem{cor}{Corollary}[section]

\newtheorem{defn}{Definition}

\def\L{\mathcal L}
\def\T{\mathcal T}

\def\CC{\mathbb C}
\def\TT{\mathbb T}
\def\ZZ{\mathbb Z}

\def\RR{\mathbb R}
\def\QQ{\mathbb Q}

\title{ Cutting Sequences and Palindromes}

\author{Jane Gilman and Linda Keen }

\address{Department of  Mathematics, Rutgers University, Newark, NJ
 07079} \email{gilman@andromeda.rutgers.edu} \thanks{Rutgers Research Council \& Yale University}
\address{Department of  Mathematics,  Lehman College and Graduate
Center,  CUNY, Bronx, 10468} \email{linda.keen@lehman.cuny.edu}

\thanks{PSC-CUNY}
\subjclass[2000]{}
\begin{document}
\maketitle

\begin{abstract}
We give a unified geometric approach to some theorems about
primitive elements and palindromes in free groups of rank 2. The
geometric treatment gives new proofs of the theorems.  This paper is
dedicated to Bill Harvey on the occasion of his 65th birthday.

\end{abstract}

\section{Introduction}

In this paper we discuss four  older more or less well-known
theorems about two generator free groups and a more recent one,
an enumerative scheme for primitive words. We describe a geometric
technique that ties all of these theorems together and gives new
proofs of four of them. This approach and the enumerative scheme
will be useful in applications. These applications will be studied
elsewhere \cite{GKgeom}.

The main object here is a two generator free group which we denote
by $G=\langle A,B \rangle$. 

\begin{defn} A word $W=W(A,B) \in G$ is {\em primitive} if there is
another word $V=V(A,B) \in G$ such that $W$ and $V$ generate $G$.
$V$ is called a {\em primitive associate} of $W$ and the unordered
pair $W$ and $V$  is called a {\em pair of primitive associates}.
\end{defn}

\begin{defn} A word $W=W(A,B) \in G$ is a {\em palindrome} if it reads the same forward and backwards.
\end{defn}

In \cite{GKwords} we found connections between a number of different
forms of primitive words and pairs of primitive associates in a two
generator free group. These were obtained using both algebra and
geometry. The theorems that we discuss, Theorems \ref{thma1},
\ref{thma2}, \ref{thma3}, \ref{thm:exponent formula}  can be found
in \cite{GKwords} and Theorem~\ref{thma4} can be found in
\cite{Piggott}, and Theorem~\ref{thma5}, the enumeration scheme,
along with another proof of Theorem~\ref{thma4} can be found in
\cite{GKenum}.

There are several different geometric objects that can be associated
 to
  two generator free groups; among them are  the punctured
torus,  the a three holed sphere and the genus two handlebody. Here
we focus on the punctured torus and use ``cutting sequences'' for
simple curves to obtain proofs of
Theorems~\ref{thma1},~\ref{thma2},~\ref{thma3} and~\ref{thma4}.

A similar treatment can be made for the three holed sphere. It was
in this setting that we first noticed that the palindromes and
products of palindromes were inherent in the geometry by looking at
the technique developed in Vidur Malik's thesis \cite{Malik} for the
three holed sphere representation of two generator groups. The
concept of a geometric center of a primitive word was inherent in
his work. We thank him for his insight.


\section{Notation and Definitions}

In this section we establish the notation and give the definitions
needed to state the five theorems and we state them. Note that in
stating these theorems in the forms below we are gathering together
results from several places into one theorem. Thus, for example, a
portion of the statements in theorem \ref{thma1} appears in
\cite{KS} while another portion appears in \cite{GKwords}.

A word $W=W(A,B) \in G$ is an expression  $A^{n_1}B^{m_1}A^{n_2}
\cdots B^{n_r}$ for some set of $2r$ integers
$n_1,...,n_r,m_1,...,n_r$.

The first theorem gives necessary conditions that
   the sequence of exponents of primitive words satisfy.
  These are called {\em primitive exponents}.
  That is, we see in Theorem~\ref{thma1} that there is a
rational number  $p/q$ that is associated to the word via its
primitive exponents.  Necessary and sufficient conditions for the
word to be primitive are given in Theorem~\ref{thm:exponent
formula}.

%
\begin{thm}\label{thma1} {\rm{(\cite{GKwords,KS})}} If $W=W(A,B)$ in $G=\langle A,B \rangle$ is primitive then up
to cyclic reduction and inverse, it has either the form

\begin{equation}\label{eqn1a} B^{n_0}A^{\epsilon}B^{n_1}A^{\epsilon}B^{n_2} \ldots
A^{\epsilon}B^{n_p}\end{equation}  where $\epsilon = \pm 1$  and  $
q=\sum_{i=1}^p n_i $ with  $p$ and $q$  relatively prime; the
exponents satisfy   $n_j =[q/p]$ or $n_j=[q/p]+1$, $ 0<j \leq p$,
where $[\;\;]$ denotes the greatest integer function,
 and no two adjacent exponents are both $[q/p]+1$;\\

 or it has the form
\begin{equation}\label{eqn2} A^{n_0}B^{\epsilon}A^{n_1}B^{\epsilon}A^{n_2} \ldots
B^{\epsilon}A^{n_q}\end{equation} where $\epsilon = \pm 1$ and    $
\sum_{i=1}^q n_i=p $ with $p$  and $q$   relatively prime; the
exponents satisfy $n_j =[p/q]$ or $n_j=[p/q]+1$, $ 0<j \leq p$,
 and no two adjacent exponents are both $[p/q]+1$.

We denote the word in either of the forms (\ref{eqn1a}) and
(\ref{eqn2}) by $W_{p/q}$. Which form is determined by whether $p/q$
is greater or less than $1$.

Two primitive words $W_{p/q}$ and $W_{r/s}$ are a pair of primitive
associates if and only if
 $|ps-qr|=1$.
\end{thm}

\subsection{Farey  arithmetic}

In what follows when we use $r/s$ to denote a rational, we assume
that $r$ and $s$ are integers, $s \ne 0$ and $(r,s)=1$.  We let
$\QQ$ denote the rational numbers, but we think of the rationals as
being points on the real axis in the complex plane.  We use the
notation $1/0$ to denote the point at infinity.

 To state
the second theorem, we need the concept of Farey addition for
fractions.

\begin{defn} If $\frac{p}{q}, \frac{r}{s} \in \QQ$, the {\em Farey sum} is
$$\frac{p}{q} \oplus \frac{r}{s} = \frac{p+r}{q+s}$$
Two fractions are called {\em Farey neighbors} or simply called {\sl
neighbors} if $|ps-qr=1|$ and the corresponding words are also
called {\sl neighbors}.
\end{defn}

 Note that the Farey neighbors
of $1/0$ are the rationals $n/1$. If $\frac{p}{q}< \frac{r}{s}$ then
it is a simple computation to see that
$$\frac{p}{q} < \frac{p}{q} \oplus \frac{r}{s} < \frac{r}{s}$$ and
that both pairs of fractions $$(\frac{p}{q},\frac{p}{q} \oplus
\frac{r}{s}) \mbox{  and  } (\frac{p}{q} \oplus \frac{r}{s},
\frac{r}{s})$$ are neighbors if $(p/q, r/s)$ are.

We can create the diagram for the Farey tree  by marking each
fraction by a point on the real line and joining each pair of
neighbors by a semi-circle orthogonal to the real line in the upper
half plane. The points $n/1$ are joined to their neighbor $1/0$ by
vertical lines. The important thing to note here   is that because
of the properties above none of the semi-circles or lines intersect
in the upper half plane. To simplify the exposition when we talk
about a point or a vertex we also mean the word corresponding to
that rational number.  Each pair of neighbors together with their
Farey sum form the vertices of a curvilinear or hyperbolic triangle
and the interiors of two such triangles are disjoint. Together the
set of these triangles forms  a tessellation of the hyperbolic plane
which is known as the Farey tree.

Let $W_{p/q}$ and $W_{r/s}$ be two primitive words labeled by
rational numbers ${\frac{p}{q}}$ and ${\frac{r}{s}}$. We can always
form the product  $W_{p/q} \cdot W_{r/s}$. If $p/q$ and $r/s$ are
neighbors, the word $W_{(p+r)/(q+s)} = W_{p/q} \cdot W_{r/s}$ so
that Farey sum corresponds to concatenation of words and by abuse of
language we talk about the Farey sum of words.


Fix  any point $\zeta$ on the positive imaginary axis. Given a
fraction, ${\frac{p}{q}}$, there is a hyperbolic geodesic $\gamma$
from $\zeta$ to ${\frac{p}{q}}$ that intersects a minimal number of
these triangles.

\begin{defn} The {\em Farey level} or the {\em level} of $p/q$, $Lev(p/q)$ is the number of triangles
traversed by $\gamma$ \end{defn}

Note that the curve (line) $\gamma$ joining $\zeta$ to either $0/1$
or $1/0$ does not cross any triangle so these rationals have level
$0$. The geodesic joining $\zeta$ to $1/1$  intersects only the
triangle with vertices $1/0, 0/1$ and $1/1$  so the level of $1/1$
is  $1$. Similarly the level of $n/1$ is $n$.

We emphasize that we now have two   different  and independent
orderings on the rational numbers: the ordering as rational numbers
and their ordering by level. That is, given ${\frac{p}{q}}$ and
${\frac{r}{s}}$, we might, for example, have as rational numbers
${\frac{p}{q}} \le {\frac{r}{s}}$, but $Lev({\frac{r}{s}}) \le
Lev({\frac{p}{q}})$. If we say one rational is larger or smaller
than the other,  we are referring to the standard  order on the
rationals. If we say one rational is higher or lower than the other,
we are referring to the {\sl levels} of the fractions.

\begin{defn} We determine a {\em Farey sequence for ${\frac{p}{q}}$} inductively  by choosing the
new vertex of the next triangle  in the sequence of triangles
traversed by $\gamma$. \end{defn}

The Farey sequence for ${\frac{3}{5}}$ is shown in Figure 1.

Given $p/q$, we can find the smallest and largest rationals $m/n$
and $r/s$ that are its neighbors.  These also have the property that
they are the only neighbors with lower level. That is, as rational
numbers  $m/n < p/q < r/s$ and the levels satisfy $Lev(m/n) <
Lev(p/q)$ and  $Lev(r/s) < Lev(p/q)$, and if $u/v$ is any other
neighbor $Lev(u/v) >   Lev(p/q)$.

\begin{defn}  We call the smallest and the largest neighbors of the rational $p/q$
 the {\em distinguished neighbors} of $p/q$.
\end{defn}

Note that we can tell whether which distinguished neighbor $r/s$ is
smaller (respectively larger) than $p/q$ by the sign of  $rq-ps$.

Farey sequences are related to continued fraction expansions of
fractions (see for example, \cite{HardyWright}).  In particular,
write  $$ {\frac{p}{q}} = 
   [a_0,\ldots,a_k]$$
  where $a_j > 0$, $j=1\ldots k$ and
for $n=0,\ldots, k-1$
  set $\frac{p_n}{q_n}=[a_0, \ldots,  a_n]$.   These
approximating fractions can be computed recursively from the
continued fraction for $p/q$ as follows:
$$p_0 = a_0,  q_0 = 1  \mbox{ and } p_1=a_0a_1+1, q_1 = a_1 $$
$$ p_j = a_j p_{j-1} + p_{j-2} \, \, , q_j = a_j q_{j-1} + q_{j-2} \, \, j=2, \dots, k.$$

The level of $p/q$ can be expressed in terms of the continued
fraction expansion  by the formula
$$Lev(p/q) = \sum_{j=0}^k a_j.$$ The distinguished neighbors of
$p/q$ have continued fractions
$$[a_0, \ldots, a_{k-1}] \mbox{  and   } [a_0, \ldots,
a_{k-1},a_k-1].$$

The Farey sequence contains the approximating fractions as a
subsequence.  The points of the Farey sequence between
$\frac{p_j}{q_j}$ and $\frac{p_{j+1}}{q_{j+1}}$ have continued
fraction expansions
$$[a_0,a_1, \ldots a_{j} + 1], [a_0,a_1, \ldots, a_{j} + 2], \ldots,  [a_0,a_1, \ldots a_{j} +
a_{j+1}-1].$$

{\em As real numbers}, the approximating fractions
$\frac{p_j}{q_j}$, termed the {\sl approximants},   are alternately
larger
 and smaller than  $\frac{p}{q}$.
 The number $a_j$ counts the number of
times the new endpoint in the Farey sequence lies on one side of the
old one.

Note that if $p/q>0$, then $0 \leq a_0 <p/q$.  The even approximants
$p_{2j}/q_{2j}$ are less than $p/q$ and the odd ones
$p_{2j+1}/q_{2j+1}$ are greater.

\subsection{Farey words, continued fraction expansions and
algorithmic words}\label{sec:2.2}

The next theorem gives a recursive enumeration scheme for primitive
words using Farey sequences of rationals.

\begin{thm}\label{thma2} {\rm{(\cite{GKwords,KS})}} The primitive words
 in $G=\langle A,B \rangle$ can be enumerated inductively  by
using Farey sequences  as follows:    set
$$W_{0/1}=A, \quad W_{1/0}=B.$$
Given $p/q$, consider its Farey sequence.  Let ${\frac{m}{n}}$ and
${\frac{r}{s}}$ be its distinguished neighbors labeled so that
$${\frac{m}{n}} < {\frac{p}{q}} < {\frac{r}{s}}.
$$ Then
$$ W_{\frac{p}{q}}=W_{\frac{m}{n} \oplus \frac{r}{s}}=W_{r/s} \cdot W_{m/n}.$$ A pair
$W_{p/q},W_{r/s}$ is a pair of  primitive associates  if and only if
$\frac{p}{q}, \frac{r}{s}$ are neighbors, that is, $|ps-qr| =1$.
\end{thm}
We use the same notation for these words as those in
Theorem~\ref{thma1} because, as we will see when we give the proofs
of the theorems, we get the same words. Since we will also introduce
two other enumeration schemes later, we will refer to this is the
$W$-enumeration scheme when clarification is needed. The other
iteration schemes will be the $V$ and the $E$-enumeration schemes.

We note that the two products   $W_{m/n} \cdot W_{r/s}$ and $W_{r/s}
\cdot W_{m/n}$ are always conjugate in $G$.   In this $W$-iteration
scheme we always choose the product where the larger index comes
first.  The point is that in order for the scheme to work the choice
has to be made consistently.  We emphasize that $W_{p/q}$ always
denotes the word obtained using this enumeration scheme.

\medskip

 The $W_{p/q}$
words can be expanded using their continued fraction exponents
instead of their primitive exponents. This is also known as the {\sl
algorithmic form}  of the primitive words, that is, the form in
which the words arise in the $PSL(2,\mathbb{R})$ discreteness
algorithm \cite{G3,Compl,YCJiang, Malik}.

The algorithm begins with a pair of generators $(X_0,Y_0)$ for a
subgroup of $PSL(2,\RR)$ and runs through a sequence of primitive
pairs of generators. At each step the algorithm replaces a
generating pair $(X,Y)$ with either $(XY,Y)$ or $(Y, XY)$ until it
arrives at a pair that stops the algorithm and prints out {\sl the
group is discrete} or {\sl the group is not discrete}. The first
type of step is termed a linear step and the second a Fibonacci
step. Associated to any implementation of the algorithm is a
sequence of integers, the $F$-sequence of Fibonacci sequence which
tells how many linear steps occur between consecutive Fibonacci
steps. The algorithm can be run backwards from the stopping
generators when the group is discrete and free and any primitive
pair can be obtained from the stopping generators using the
backwards $F$-sequence. The $F$-sequence of the algorithm is used in
\cite{G3} and \cite{YCJiang} to determine the computational
complexity of the algorithm. In \cite{G3} it is shown that most form
of the algorithm are polynomial and in \cite{YCJiang} it is shown
that all forms are.

In \cite{GKwords} it is shown that the $F$-sequence that determines
a primitive word is equivalent to the continued fraction expansion
of the rational corresponding to that primitive word. The following
theorem exhibits the primitive words with the continued fraction
expansion exponents in its most concise form.

\begin{thm}\label{thma3} {\rm{(\cite{GKwords})} } If $[a_0, \ldots, a_k]$ is the continued fraction expansion of $p/q$,
the primitive word $W_{p/q}$ can be written inductively using the
continued fraction  approximants be $p_j/q_j=[a_0, \ldots, a_j]$.
They are alternately larger and smaller than $p/q$.

Set
 $$W_{0/1}=A, \, \,   W_{1/0}=B  \mbox{  and   } W_{1/1}=BA.$$ For
 $j=1, \ldots, k$ if
$p_{j-2}/q_{j-2}> p/q$  set
$$W_{p_j/q_j}=W_{p_{j-2}/q_{j-2}}(W_{p_{j-1}/q_{j-1}})^{a_j}$$ and set
$$W_{p_j/q_j}=(W_{p_{j-1}/q_{j-1}})^{a_j}W_{p_{j-2}/q_{j-2}}$$
otherwise.
\end{thm}

We have an alternative recursion which gives us formulas for the
primitive exponents and hence necessary and sufficient conditions to
recognize primitive words.

Assume $p/q>1$ and write $p/q=[a_0, \ldots a_k]$. By assumption
$a_0>0$.
   If $0<p/q<1$ interchange $A$ and $B$ and argue
  with $q/p$.

Set $V_{-1}=B$ and $V_0=W_{p_0/q_0}=AB^{a_0}$. Then for $j=1,
\ldots, k$, $V_j=V_{j-2}[V_{j-1}]^{a_j}$.

\begin{thm}\label{thm:exponent formula} {\rm{(\cite{GKwords, Malik})}}
 Write
$$V_j=   B^{n_0(j)}AB^{n_1(j)}\ldots AB^{n_{t_j}(j)}$$ for $j=0, \ldots,
k$.  The primitive exponents of $W_{p/q}$ are related to the
continued
fraction of $p/q=[a_0, \ldots, a_k]$ as follows:\\

If $j=0$, then $t_0=1=q_0$,  $n_0=0$ and  $n_1(0)=a_0$.
If $j=1$, then $t_1=a_1=q_1$, $n_0(1)=1$ and $n_i(1)=a_0$, $i=1,
\ldots, a_1$.


If $j=2$ then $t_2=a_2a_1+1=q_2$,    $n_0(2)=0$ and $n_i(2)=a_0+1$
for $i \equiv \, 1 \, \mbox{ mod } \, q_2$ and $n_i(j) = a_0$
otherwise.

For $j>2, \ldots, k$,
\begin{itemize}
\item   $n_0(j)=0$ if $j$ is even and $n_0(j)=1$ if $j$ is odd.
\item $t_j=q_j$.
\item  For $i=0 \ldots t_{j-2}$, $n_i(j)=n_i(j-2)$.
\item  For $i=t_{j-2}+1 \ldots t_j$,  $n_i(j)=a_0+1$ if $i \equiv
t_{k-2}  \mbox{ mod } t_{k-1}$ and $n_i(j) = a_0$ otherwise.
\end{itemize}
These conditions on the exponents are necessary and sufficient for a
word (up to cyclic permutation) to be primitive.
\end{thm}

Using the recursion formulas, we obtain a new proof of the following
corollary which was originally proved in
 \cite{BuserS}. We omit the proof as it is a fairly straightforward induction argument on the Farey level.

  \begin{cor}  {\rm({\cite{BuserS})}} In the expression for a primitive word $W_{p/q}$, for any integer
$m$, $0<m<p$, the sums of any $m$ consecutive primitive exponents
$n_i$
  differ by at most $\pm 1$.  \end{cor}

The following theorem was originally  proved in \cite{OZ} and in
\cite{Piggott} and \cite{KR}.

\begin{thm}\label{thma4} {\rm {(\cite{GKenum,KR,OZ,Piggott})}} Let $G=\langle A,B  \rangle$ be a two generator free
group.  Then any primitive element $W \in G$ is conjugate to a
cyclic permutation of either a palindrome in $A,B$ or a product of
two palindromes. In particular, if the length of $W$ is $p+q$, then,
up to cyclic permutation, $W$ is a palindrome if and only if $p+q$
is odd and is a product of two palindromes otherwise.
\end{thm}

We note that this can be formulated equivalently using the parity of
$pq$ which is what we do below.

  In the $pq$ odd case, the two palindromes in the
product can be chosen in various ways. We will make a particular
choice in the next theorem.

\subsection{$E$-Enumeration}

The next theorem, proved in \cite{GKenum},  gives yet another
enumeration scheme for primitive words, again using Farey sequences.
The new scheme to enumerate primitive elements is useful in
applications, especially geometric applications. These applications
will be studied elsewhere \cite{GKgeom}.  Because the words we
obtain are cyclic permutations of the words $W_{p/q}$, we use a
different notation for them; we denote them as $E_{p/q}$.

\begin{thm}\label{thma5}  {\rm{(\cite{GKenum})}} The primitive elements of a two generator
free group can be enumerated recursively using their Farey sequences
as follows.  Set $$E_{0/1} = A, \, \, E_{1/0}=B, \,\, \mbox{  and  }
E_{1/1}=BA.$$

Given $p/q$ with distinguished neighbors $m/n,r/s$ such that
$m/n<r/s$, \begin{itemize}\item if $pq$ is odd, set
$E_{p/q}=E_{r/s}E_{m/n}$ and
\item  if $pq$ is even, set
$E_{p/q}=E_{m/n}E_{r/s}$.  In this case $E_{p/q}$ is the unique
palindrome cyclicly conjugate to $W_{p/q}$.
\end{itemize}

We also use $P_{p/q}$ for $E_{p/q}$ when $pq$ is even and $Q_{p/q}$
when $pq$ is odd.

 $E_{p/q}$ and $E_{p'/q'}$ are primitive
associates if and only if $pq'-qp'=\pm 1$.
\end{thm}

Note that when $pq$ is odd, the order of multiplication is the same
as in the enumeration scheme for $W_{p/q}$ but when $pq$ is even, it
is reversed. This theorem says that if $pq$ is even, $E_{p/q}$ is
the  unique palindrome cyclicly conjugate to $W_{p/q}$.  If $pq$ is
odd, then $E_{p/q}$  determines a canonical factorization of (the
conjugacy class of) $W_{p/q}$ into a pair of palindromes. This
factorization exhibits the Farey sequence of $p/q$ and the order of
multiplication is what makes the enumeration scheme work.

In this new enumeration scheme, Farey neighbors again correspond to
primitive pairs but the elements of the pair  $(W_{p/q},W_{p'/q'})$
are not necessarily conjugate to the elements of the pair $(E_{p/q},
E_{p'/q'})$ by the same element of the group. That is, they are not
necessarily conjugate as pairs.

 \section{Cutting Sequences}

 We represent $G$ as the fundamental group of a punctured torus
and use the technique of {\sl cutting sequences} developed by Series
(see \cite{S, KS, Nielsen}) as the unifying theme. This
representation assumes that the group $G$ is a discrete free group.
Cutting sequences are a variant on Nielsen boundary development
sequences \cite{Nielsen}. In this section we outline the steps to
define cutting sequences.

\begin{itemize}
\item It is  standard that $G=\langle
A,B  \rangle$  is isomorphic to the fundamental group of a punctured
torus $S$.   Each element of $G$ corresponds to a free homotopy
class of curves on $S$.  The primitive elements are classes of
simple curves that do not simply go around the puncture.  Primitive
pairs are classes of simple closed curves with a single intersection
point.
\item Let $\L$ be the lattice of points in $\CC$ of the form $m+ni, \, m,n \in
\ZZ$ and let $\T$ be the corresponding lattice group generated by
$a=z \mapsto z+1, b=z \mapsto z+i$.  The (unpunctured)  torus is
homeomorphic to the quotient  $ \TT= \CC   / \T$.  The horizontal
lines map to longitudes  and the vertical lines to meridians on
$\TT$.

  The punctured torus is homeomorphic to the quotient of the
plane punctured at the lattice, $(\CC \setminus \L) / \T$.  Any
  curve in $\CC$ whose endpoints are identified by the commutator
  $aba^{-1}b^{-1}$ goes around a puncture and is no longer
  homotopically trivial.

\item  The simple closed curves on $\TT$ are exactly the projections of lines
joining pairs of lattice points (or lines parallel to them). These
are lines $L_{q/p}$ of rational slope $q/p$.
 The projection $l_{q/p}$ consists of $p$ longitudinal
loops   $q$ meridional loops.  We assume that  $p$ and $q$ are
relatively prime; otherwise the curve has multiplicity equal to the
common factor.

  For the punctured torus, any line of rational slope, not passing through the punctures
  projects to a simple closed curve and any simple closed curve, not enclosing the puncture, lifts to a curve
  freely homotopic to a line of rational slope.

\item Note that, in either case, if we try   to draw
  the projection of $L_{q/p}$ as a simple curve, the order in which we traverse the
  loops on $\TT$ (or $S$) matters. In fact there is, up to cyclic permutation and reversal, only
  one way to draw the curve.  We will find this way using cutting
  sequences. Below, we assume we are working on $\TT$.

\item  Choose as fundamental domain (for $S$ or $\TT$)  the square $D$ with corners (puncture
points) $\{0,1,1+i,i\}$.  Label the bottom side $B$ and the top side
$\bar{B}$; label the left side $A$ and the right side $\bar{A}$.
Note that the transformation $a$ identifies  $A$ with $\bar{A}$ and
$b$ identifies $B$ with $\bar{B}$.  . Use the translation group to
label the sides of all copies of $D$ in the plane.

\item
Choose a fundamental segment of the line $L_{q/p}$ and pick one of
its endpoints as initial point. It passes through $p+|q|$ copies of
the fundamental domain. Call the segment in each copy a strand.

Because the curve is simple, there will either be ``vertical''
strands joining the sides $B$ and $\bar{B}$, or ``horizontal''
strands joining the sides $A$ and $\bar{A}$, but not both.

 Call the segments joining a horizontal and vertical side corner
strands.  There are four possible types of corner strands: from left
to bottom, from left to top, from bottom to right, from top to
right. If all four types were to occur, the projected curve would be
trivial on $\TT$. There cannot be only one or three different types
of corner strands because the curve would not close up. Therefore
the  only corner strands occur on one pair of opposite corners and
there are an equal number on each corner.

\item Traversing  the fundamental segment from its initial point, the
line goes through or ``cuts'' sides of copies of $D$.  We will use
the side labeling to define a {\em cutting sequence} for the
segment. Since each side belongs to two copies it has two labels. We
have to pick one of these labels in a consistent way. As the segment
passes through, there is the label from the copy it leaves and the
label from the copy it enters. We always choose the label from the
copy it enters.  Note that the cyclic permutation depends on the
starting point.

\item
 If $|q|/p <1$, the resulting cutting sequence will contain $p$
$B$'s (or $p$ $\bar{B}$'s), $|q|$ $B$'s (or $|q|$  $\bar{B}$'s) and
there will be $p-|q|$ horizontal strands and $p$ corner strands; if
$|q|/p >1$, the resulting cutting sequence will contain $p$ $B$'s
(or $p$ $\bar{B}$'s), $|q|$ $A$'s (or $|q|$  $\bar{A}$'s) and there
will be $ |q| -p$ vertical strands and $|q|$ corner strands. We
identify the cutting sequence with the word in $W$ interpreting the
labels $A,B$ and the generators and the labels $\bar{A},\bar{B}$ as
their inverses.

\item Given an arbitrary  word $W=A^{m_1}B^{n_1}A^{m_2}B^{n_2} \ldots
A^{m_p}B^{n_p}$ in $G$, we can form a cutting sequence for it  by
drawing
  strands from the word through successive copies of $D$.  Consider
  translates of the resulting curve by elements of the lattice
  group.  If they are all disjoint up to homotopy, the word is
  primitive.
\end{itemize}

Let us illustrate with three examples. In the first two, we draw the
cutting sequences for the fractions $\frac{q}{p}=\frac{1}{1}$ and
$\frac{q}{p}=\frac{3}{2}$. In the third, we construct the cutting
sequence for the word $A^2B^3$.
\begin{itemize}

\item   A fundamental segment of $l_{1/1}$ can be chosen to begin at a point on the left ($A$) side
and pass through  $D$ and the adjacent copy above $D$;
    There will be a single corner   strand connecting  the $A$ side
    to  a
  $\bar{B}$ side and 
    another   connecting a  $B$ side to an  $\bar{A}$ side.

  To read off the   cutting sequence    begin with the point on
  $A$ and write $A$.  Then as we enter the next (and last) copy of
  $D$ we have an $B$ side.  The word is thus $ AB$.

  Had we started on the bottom, we would have obtained the word
  $BA$.

\item
A fundamental segment of $L_{3/2}$ passes through $5$ copies of the
fundamental domain. (See Figure 2.) 
There is one
 ``vertical'' segment joining  a $B$ and a $\bar{B}$, 2 corner segments joining an $A$ and a $\bar{B}$
  and two joining
   the opposite corners.  Start on  the left side.  Then, depending
   on where on this side we begin we obtain the word $ABABB$ or
   $ABBAB$.

   If we start on the bottom so that the vertical side is in the
   last copy we encounter we get $BABAB$.

\item  To see that the word $AABBB$ cannot correspond to a simple
loop, draw the a vertical line of length $3$ and join it to a
horizontal line of length $2$.  Translate it one to the right and
one up.  Clearly the translate intersects the curve and projects to
a self-intersection on the torus.   This will happen whenever the
horizontal segments are not separated by a vertical segment.

Another way to see this is to try to draw  a curve with  3 meridian
loops and two longitudinal loops on the torus. You will easily find
that if you try to connect them arbitrarily the strands will cross
on $\TT$, but if you use the order  given by the cutting sequence
they will not. Start in the middle of the single vertical strand and
enter a letter every time you come to the beginning of a new strand.
We get $BABAB$.

\item Suppose $W= B^3A^2$.  To draw the cutting sequence, begin on
the bottom of the square and, since the next letter is $B$ again,
draw a vertical strand to a point on the top and a bit to the right.
Next, since we have a third $B$, in the copy above $D$ draw another
vertical strand to the top and again go a bit to the right. Now the
fourth letter is a $A$ so we draw a corner strand to the right.
Since we have another $A$ we need to draw a horizontal strand.  We
close up the curve with a last corner strand from the left to the
top.

Because we have both horizontal and vertical strands, the curve is
not simple and the word is not primitive.

\end{itemize}

\section{Proofs}

\medskip \noindent{\bf Proof of Theorem~\ref{thma1} and Theorem~\ref{thm:exponent formula}.}
We have seen that a word is primitive if and only if its cutting
sequence has no intersecting strands and corresponds to a line of
rational slope $q/p$. We want to examine what the cutting sequences
look like for these lines.

   The cases $p/q=0/1, 1/0$ are trivial. Suppose first that $q/p \geq 1$.
 The other cases follow in the
same way, either interchanging $A$ and $B$ or replacing $B$ by
$\overline{B}$.

The line $L_{q/p}$ has slope at least $1$ so there will be at least
one vertical strand and
  no horizontal strands. Set $q/p=[a_0, \ldots a_k]$. Since $q/p>1$ we know that  $a_0>0$.

Note that there is an
    ambiguity in this
representation; we have $[a_0, \ldots
  a_k-1,1]=[a_0, \ldots a_k]$.  We  can eliminate this by
  assuming $a_k>1$.  With this convention, the parity of $k$ is well defined.

  Assume first $k$ is even, choose as starting point the
lowest point on an $A$ side.   Because there are no horizontal
strands, we must either go up or down; assume we go up. The first
letter in the cutting sequence is $A$ and since the strand must be a
corner strand, the next letter is $B$. As we form the cutting
sequence we see that because there are no horizontal strands, no $A$
can be followed by another $A$. Because we started at the lowest
point on $A$, the last strand we encounter before we close up must
start at the rightmost point on a $B$ side. Since there are $p+q$
strands, this means the sequence, and hence the word has the first
form of Theorem~\ref{thma1}. Since $p/q>0,$   in the exponents of
the $A$ terms $\epsilon=1$.  Since we begin with an $A$, $n_0=0$ and
$$W_{p/q}= AB^{n_1}AB^{n_2}A \ldots B^{n_p}, \, \sum n_i = q. $$
If we use the translation group to put all the strands into one
fundamental domain, the endpoints of the strands on the sides are
ordered. We see that if we are at a point on the $B$ side, the next
time we come to the  $B$ side we are at a point that is $p$ to the
right $\rm{ mod}(q)$.

Let us see exactly what the exponents are. Since we began with the
lowest point on the left, the first $B$ comes from the $p^{th}$
strand on the bottom.  There are $q$ strands on the bottom; the
first (leftmost) $q-p$ strands are vertical and   the last $p$ are
corner strands.  Since we move to the right $p$ strands at a time,
we can do this $a_0=[q/p]$ times. The word so far is $AB^{a_0}$.

 At this point we have a corner strand so the next letter will be an $A$.
Define $r_1$ by $q=a_0 p + r_1$.  The corner strand ends at the
right endpoint $r_1+1$ from the bottom and the corresponding corner
strand on the $A$ side joins with   the $(p-r_1)^{th}$ vertical
strand on the bottom.  We again move to the right $p$ strands at a
time, $a_0$ times, while $a_0p-r_1> q-p$. After some number of
times, $a_0p-r_1 \leq q-p$. This number, $n$,  will satisfy $p= r_1
n +r_2$ and $r_2<r_1$. Notice that this is the first step of the
Euclidean algorithm for the greatest common denominator and it
generates the continued fraction coefficients at each step. Thus
$n=a_1$ and  the word at this point is $[AB^{a_0}]^{a_1}$. Since we
are now at a corner strand, the next letter is an ``extra'' $B$.  We
repeat the sequence we have already obtained   $a_3$ times where
$r_1= a_3 r_2 +r_3$ and $r_3<r_2$. The word at this point is
$[AB^{a_0}]^{a_1}B[[AB^{a_0}]^{a_1}]^{a_3}$ which is the word we
called $V_3$ in Theorem~\ref{thm:exponent formula}.

We continue in this way.    We see that the Euclidean algorithm
tells us that each time we have an extra $B$ the sequence up to that
point repeats  as many times as the next $a_i$  entry in the
continued fraction expansion of $q/p$.  When we come to the last
entry $a_k$, we have used all the strands and are back to our
starting point. We see that  the exponent structure is forced on us
by the number $q/p$ and the condition that the strands not
intersect.

If $k$ is odd, we begin the process at the rightmost bottom strand
and begin the word with $B$ and obtain the recursion.

 Note that had
we chosen a different starting point we would have obtained a cyclic
permutation of $W_{q/p}$, or, depending on the direction, its
inverse.

Thus, if the exponents $n_i$ of a word $W$ with $\sum_{i=1}^p n_i
=q$, or some cyclic permutation of it, do not satisfy these
conditions, the strands of its cutting sequence must either
intersect somewhere or they do not close up properly and the word is
not primitive. The conditions are therefore both necessary and
sufficient for the word to be primitive and
Theorem~\ref{thm:exponent formula} follows.

It is obvious that the only primitive exponents that can occur are
$a_0$ and $a_0+1$.  Moreover, no adjacent primitive exponents can
equal $a_0+1$.  This gives the simple necessary conditions of
Theorem~\ref{thma1}.

The  primitive exponent formulas in Theorem~\ref{thm:exponent
formula} follow by induction on $k$.

  For  $0<q/p<1$ we have no vertical strands and we interchange the
 roles of $A$ and $B$.  We use the continued fraction
 $p/q=[a_0,\ldots,a_k]$ and argue as above, replacing ``vertical''
 by ``horizontal''.

  For $p/q<0$, we replace $A$ or $B$ by
 $\bar{A}$ or $\bar{B}$ as appropriate.

 To see when two primitive words $W_{p/q}$ and $W_{r/s}$ are associates,
 note that the lattice $\L$ is generated by
fundamental segments of lines $L_{p/q}, L_{r/s}$ if and only if
$|ps-qr|=1$, or equivalently, if and only if $(p/q,r/s)$ are
neighbors. \hfill $\Box$

 \noindent{\bf Proof of   Theorem~\ref{thma2} and \ref{thma3}.}
Although Theorem~\ref{thma2} and \ref{thma3} can be deduced from the
proof above, we give an independent proof.

The theorems prescribe  a recursive definition of a primitive word
associated to a rational $p/q$. We assume $m/n$ and $r/s$ are
distinguished neighbors and $$\frac{m}{n} < \frac{p}{q} <
\frac{r}{s}.$$

 We need to show that if   we draw the strands for the cutting sequence
for $(W_{m/n}$ and $W_{r/s})$  in the same diagram, then the result
is the cutting sequence of the product.

Note first that if $r/s,m/n$ are   neighbors, the vectors joining
zero  with $m+ni$ and $r+si$
    generate the lattice $\L$.  Draw a fundamental segment $s_{m/n}$ for
    $W_{m/n}$ joining $0$ to $m+ni$ and a fundamental segment $s_{r/s}$ for $W_{r/s}$ joining
    $m+ni$ to $(m+r)+(n+s)i$.  The straight line $s$ joining $0$ to
    $(m+r)+(n+s)i$ doesn't pass through any of the lattice points because by
    the neighbor condition $rn-sm=1$, $s_{m/n}$ and $s_{r/s}$ generate the lattice.  We   therefore get the same
    cutting sequence whether we follow  $s_{m/n}$ and $s_{r/s}$  in turn or follow the straight line $s$.
    This means that the cutting sequence for $W_{p/q}$ is the
    concatenation of the cutting sequences of $W_{r/s}$ and
    $W_{m/n}$  which is what we had to show.

    This observation about the generators of the lattice also proves
    that if $r/s,m/n$ are
neighbors, the pair $W_{r/s},W_{m/n}$ is a pair of primitive
associates.

We note that proving Theorem~\ref{thma3} is just a matter of
notation.

\hfill $\Box$

Notice that this theorem says that,  if given a   primitive
associate pair $(W_{p/q},W_{r/s})$,  we draw the strands for cutting
sequence for each primitive in the same diagram, then the  result is
the cutting sequence of the product.

\medskip \noindent {\bf Proof of Theorem~\ref{thma4}.}

Suppose  $pq$ is even.  Again we prove the theorem for $0<p/q<1$.
The other cases follow as above by interchanging the roles of $A$
and $B$ or replacing $B$ by $\overline{B}$.  The idea is to choose
the starting point correctly.

Draw  a line of slope $p/q$.  By assumption, there are horizontal
but no vertical strands and $p-q>0$ must be odd. This implies that
in a fundamental segment there are an odd number of horizontal
strands. In particular, if we pull all the strands of a fundamental
segment into one copy of $D$, one of the horizontal strands is the
middle strand. Choose the fundamental segment for the line in the
lattice so that it is centered about
 this middle horizontal strand.

 To form the cutting sequence for the corresponding  word $W $,
  begin at the right endpoint of the middle strand and
 take as initial point the leftpoint that it corresponds to.  Now go to the other end of the
 middle strand on the left and take as initial point the rightpoint that it corresponds to form  the cutting sequence for a
 word $V$.  By the symmetry, since we began with a middle strand, $V $ is $W $ with all the $A$'s replaced
by $\overline{A}$'s and all the $B$'s replaced by $\overline{B}$'s.
Since $V=W ^{-1}$, we see that $W$ must be a palindrome which we
denote as    $W=P_{p/q}$. Moreover, since it is the cutting sequence
of a fundamental segment of the line of slope $p/q$,
 it must is a cyclic permutation of $W_{p/q}$.

 Note that since we began with a horizontal strand, the
first letter in the sequence is an $A$ and, since it is a
palindrome, so is the last letter.

When $q/p>1$, there are horizontal and no vertical strands, and
there is a middle horizontal strand.  This time we choose this
strand and go right and left to see that we get a palindrome. The
first and last letters in this palindrome will be $B$.

If $p/q<0$, we argue as above but either $A$ or $B$ is replaced by
respectively $\bar{A}$ or $\bar{B}$. \hfill  $\Box$

We now turn to the enumeration scheme:

\vskip .2in

\noindent {\bf Enumeration for Theorem~\ref{thma5}.}

The proof of the enumeration theorem involves purely algebraic
manipulations and can be found in \cite{GKenum}. We do not reproduce
it here but rather give a heuristic geometric  idea of the
enumeration and the connection with palindromes that comes from the
$PSL(2, \mathbb{R})$ discreteness algorithm \cite{G2,G3}.

Note that the absolute value of the trace of an element $X \in
PSL(2, \mathbb{R})$, $|\mbox{trace}(X)|$,  is well-defined. Recall
that $X$ is elliptic if $|\mbox{trace}(X)| < 2$ and hyperbolic if
$|\mbox{trace}(X)| > 2$.     As an isometry of the upper half plane,
each hyperbolic element has an invariant geodesic called its {\em
axis}. Each point on the axis is moved a distance $l(X)$ towards one
endpoint on the boundary.  This endpoint is called the attractor and
the distance can be computed from the trace by the formula $\cosh
{\frac{l(X)}{2}} = {\frac{1}{2}} |\mbox{trace}(X)|$.  The other
endpoint of the axis is a repeller.

For convenience we use the unit disk model and consider elements of
$PSL(2, \mathbb{R})$ as isometries of the unit disk. In the
algorithm one begins with a representation of the group where the
generators $A$ and $B$ are (appropriately ordered) hyperbolic
isometries of the unit disk.  The algorithm applies to any
non-elementary representation of the group where the representation
is not assumed to be free or discrete. The axes of $A$ and $B$ may
be disjoint or intersect. We illustrate the geometric idea using
intersecting axes.

If the axes of $A$ and $B$ intersect, they intersect in unique point
$p$. In this case one does not need an algorithm to determine
discreteness or non-discreteness as long as the multiplicative
commutator, $ABA^{-1}B^{-1}$,  is not an elliptic isometry.
 However, the geometric steps used in
determining discreteness or non-discreteness in the case of an
elliptic commutator still make sense. We think of the representation
as being that of a punctured torus group when the group is discrete
and free.

Normalize at the outset so that the translation
length of $A$ is smaller than the translation length of $B$,
the axis of $A$ is the geodesic joining $-1$ and $1$ with attracting fixed
point $1$ and the axis of $B$ is the line joining $  e^{i\theta}$ and $  -e^{i\theta}$. This makes  the point $p$   the origin.  Replacing $B$ by its inverse if necessary, we may assume the attracting fixed point of $B$ is $e^{i\theta}$ and $-\pi/2<\theta  \leq \pi/2$.

{\sl The geometric property of the palindromic words is that their
axes all pass through the origin.}

Suppose  $(p/q,p'/q')$ is a pair of neighbors with $pq$ and $p'q'$
even and $p/q<p'/q'$.  The word $W_{r/s}=W_{p'/q'}W_{p/q}$  is not a
palindrome or conjugate to a palindrome. Since it is a primitive
associate of both $W_{p'/q'}$ and $W_{p/q}$ the axis of
$Ax_{W_{r/s}}$ intersects each of the axes $Ax_{W_{p/q}}$ and
$Ax_{W_{p'/q'}}$ in a unique point; denote these points by $q_{p/q}$
and  $q_{p'/q'}$ respectively.  Thus, to each triangle,
$(p/q,r/s,p'/q')$ we obtain a triangle in the disk with vertices
$(0,q_{p/q}, q_{p'/q'})$.

The algorithm provides a method of choosing a next neighbor and next
associate primitive pair so that  at each step the longest side of
the triangle is replaced by a shorter side.  The procedure stops
when the sides are as short as possible. Of course, it requires
proof to see that this procedure will stop and thus will actually
give an algorithm.

There is a similar geometric description of the algorithm and
palindromes in the case of disjoint axes.

\bibliographystyle{amsplain}

\end{document}